\documentclass[12pt]{article}
\usepackage[dvips]{graphics,color}

\usepackage{latexsym}
\usepackage{amssymb}
\usepackage{amsmath}
\let\cal\mathcal

\usepackage{latexsym}
\usepackage{amssymb}
\usepackage{euscript}
\usepackage[dvips]{graphics}
\usepackage{epsf}

\let\cal=\mathcal      

\def\mcc{M\raise.5ex\hbox{c}C}
\def\mccarthy{M\raise.5ex\hbox{c}Carthy}
\def\Hu{\H}
 
\def\szs{Szeg\Hu{o}'s }

\def\ie{{\it i.e. }}

\def\atn{{\A}^2(\nu)}

\def\ptm{P^2(\mu)}
\def\ltm{L^2(\mu)}
\def\h{{\cal H}}

\def\K{{\cal K}}

\def\l{\lambda}
\def\z{\zeta}

\let\i=\infty

\def\={\ = \ }

\def\E{E_\l}

\def\C{\mathbb C}
\def\R{\mathbb R}
\def\T{\mathbb T}
\def\D{\mathbb D}

\def\Z{\mathbb Z}

\def\inn{\ \in \ }

\def\dis{\displaystyle}

\def\be{\setcounter{equation}{\value{theorem}} \begin{equation}}
\def\ee{\end{equation} \addtocounter{theorem}{1}}
\def\beq{\begin{eqnarray*}}
\def\eeq{\end{eqnarray*}}
 
\def\att{\addtocounter{theorem}{1}}
\def\vs{\vskip 5pt}
\def\bs{\vskip 12pt}

\def\exam{\bs \att {{\bf Example \thetheorem \ }} }

\def\rem{\bs \att {{\bf Remark \thetheorem \ }} }
\def\bp{{\sc Proof: }}
\def\ep{{}{\hfill $\Box$} \vskip 5pt \par}

\def\oec{{}{\hspace*{\fill} $\lhd$} \vskip 5pt \par}
\def\bl{\begin{lemma}}
\def\el{\end{lemma}}
\def\bt{\begin{theorem}}
\def\et{\end{theorem}}
\def\bprop{\begin{prop}}
\def\eprop{\end{prop}}
\def\bd{\begin{definition}}
\def\ed{\end{definition}}
\def\br{\begin{remark}}
\def\er{\end{remark}}
\def\bexer{\begin{exercise}}
\def\eexer{\end{exercise}}
\def\bfig{\begin{figure}}
\def\efig{\end{figure}}

\newtheorem{theorem}{Theorem}[section]
\newtheorem{prop}[theorem]{Proposition}
\newtheorem{lemma}[theorem]{Lemma}
\newtheorem{cor}[theorem]{Corollary}

\newtheorem{question}[theorem]{Question}
\newtheorem{definition}[theorem]{Definition}

\renewcommand{\z}{\zeta}

\renewcommand{\L}{{\mathcal L}}
\renewcommand{\E}{{\mathbb E}}
\newcommand{\mult}{\rm{mult}}
\renewcommand{\O}{\Omega}
\renewcommand{\atn}{A^2(\nu)}
\newcommand{\ath}{A^2_h(\nu)}
\newcommand{\athn}{A^2_h(\nu)}

\newcommand{\ato}{A^2(\omega)}
\newcommand{\ats}{A^2(\sigma)}
\newcommand{\om}{\omega}
\newcommand{\ran}{{\rm ran\,}}
\renewcommand{\Im}{{\rm Im\,}}
\newcommand{\ip}[2]{\langle #1, #2 \rangle}

\begin{document}
\setlength{\baselineskip}{21pt}
\title{Algebraic pairs of isometries}
\author{Jim Agler
\thanks{Partially supported by National Science Foundation Grant
DMS 0400826}\\
U.C. San Diego\\
La Jolla, CA 92093
\and
Greg Knese\\
U.C. Irvine\\
Irvine CA 92697
\and
John E. M\raise.5ex\hbox{c}Carthy
\thanks{Partially supported by National Science Foundation Grant
DMS 0501079}\\
Washington University\\
St. Louis, MO 63130}

\bibliographystyle{plain}
\maketitle
\begin{abstract}
We consider pairs of commuting isometries that are annihilated by a polynomial.
We show that the polynomial must be inner toral, which is a geometric condition on
its zero set. We show that cyclic pairs of commuting isometries are nearly unitarily
equivalent if they are annihilated by the same minimal polynomial.
\end{abstract}

\baselineskip = 18pt

\setcounter{section}{-1}
\section{Introduction}\label{seca}

Isometries form one of the best-understood classes of operators on
Hilbert spaces.  By the von Neumann-Wold decomposition, every isometry
is the direct sum of a unitary operator and a vector-valued 
shift. The non-unitary part of the isometry is called the {\em pure}
part.

Pairs of commuting isometries are more complicated.  
%
If the first isometry is pure, it can be modeled as a 
vector-valued shift, multiplication by the coordinate function on $H^2 \otimes \L$, where
$\L$ is a Hilbert space of the appropriate dimension, and $H^2$ is the Hardy space.
The second isometry then becomes multiplication
by an operator-valued inner function on $\L$, \ie 
an analytic operator-valued function on the unit disk $\D$ whose boundary values are isometric a.e.
\cite{hel64, szn-foi}. 

Although this description is very powerful, it leaves open many
questions.  The purpose of this note is to study a restricted class of
pairs of commuting isometries $V=(V_1,V_2)$, namely ones that satisfy
an algebraic relation: $q(V) = 0$ for some polynomial $q$ of two
variables.  We shall call such a pair an {\em algebraic isopair}, and
we shall say that an isopair is {\em pure} if both isometries are
pure.
Pure algebraic isopairs turn out to have a rich structure.

It is 
easy to find an algebraic isopair annihilated by the polynomial
$z^2-w^2$, but a moment's thought shows that none can be annihilated by
$z^2 - 2w^2$. 
The polynomial $1-zw$ can annihilate an isopair, but only if this is a pair
of unitaries whose joint spectrum is contained in 
$$ \T^2 \cap \{ (z,w) \ : \ 1 -zw \, = \, 0 \}. $$
(Throughout the paper, we shall use the notation that
 $\D$ is the open unit disk $\{ z \, : \, |z| < 1 \}$, $\T$ is the unit circle 
$\{ z \, : \, |z| = 1 \}$, 
and $\E$ is the exterior of the closed disk 
$\{ z \, : \, |z| > 1 \}$.)
No pure isopair is annihilated by $1-zw$.

What polynomials $q$ can 
be the minimal annihilating polynomial for some  pure isopair?

{\bf Theorem~\ref{thmb4}:} {\em Let $V=(V_1,V_2)$ be a pure algebraic isopair on
  a Hilbert space $\h$. Then there exists a square-free inner toral
  polynomial $q$ that annihilates $V$.  Moreover, if $p$ is any
  polynomial that annihilates $V$, then $q$ divides $p$.  }

A polynomial $q$ 
is called an {\em inner toral polynomial} if its 
zero set lies in $\D^2 \cup \T^2 \cup \E^2$; 
the zero set of an inner toral polynomial is called
a {\em distinguished variety}.
We discuss these in Section~\ref{secb}. 

Theorem~\ref{thmb4} gives a way to construct algebraic isopairs.
Start with an inner toral polynomial $q$; put a nice measure $\mu$ on
$Z_q \cap \T^2$; construct the Hardy space $H^2(\mu)$ that is the closure in $L^2(\mu)$
of the polynomials; and look at the pair of operators on $H^2(\mu)$ given by multiplication 
by the coordinate functions. In a way that will be made precise in Section~\ref{secc}, 
this construction in some sense gives you all cyclic algebraic isopairs.

However, they also arise in another setting.
In~\cite{fed91,veg07}, it is 
shown that on every finitely connected planar domain $R$ there
is a pair of inner functions $ (u_1,u_2)$ that map the domain conformally onto
some distinguished variety intersected with the bidisk. 
If $\nu$ is a measure
on $\partial R$ that is a log-integrable weight times harmonic measure, one 
can form a Hardy space $H^2(\nu)$
(provided every component in the complement of $R$ has interior,
this is just the closure in $L^2(\nu)$ of all functions analytic in a neighborhood of $\overline{R}$).
Multiplication by $u_1$ and $u_2$ on $H^2(\nu)$ then give a pure cyclic algebraic isopair.

\vs

In Section~\ref{secba0}, we show that a $q$-isopair (an isopair
annihilated by $q \in \C[z,w]$) can almost be broken up into a direct
sum of isopairs corresponding to each of the irreducible factors of
$q$. Specifically, we have:

{\bf Theorem~\ref{thmpropb5}:}
{\em
Let $V = (V_1,V_2)$ be a pure algebraic isopair with
  minimal polynomial $q$, and let $q_1,q_2, \dots, q_N$ be the
  (distinct) irreducible factors of $q$.  If both $V_1$ and $V_2$ have finite
  dimensional cokernels, then $V$ has a finite codimension invariant subspace
  $\mathcal{K}$ on which
\[
V\mid_{\mathcal{K}} = W_1 \oplus W_2 \oplus \cdots \oplus W_N
\]
where $W_j$ is a $q_j$-isopair, $j=1,\dots, N$.
}

The restriction to $\mathcal{K}$ is essential. Our main result says that
any two pure cyclic algebraic isopairs are nearly unitarily equivalent if and only if they
have the same minimal polynomial. ``Nearly'' means after restricting to a finite codimensional 
invariant subspace. So we say that two pairs are nearly unitarily equivalent if and only if each one is
unitarily equivalent to the other restricted to a finite codimensional 
invariant subspace. We say a pair is nearly cyclic if, when restricted to a finite codimensional 
invariant subspace, it becomes cyclic. We have:

{\bf Theorem~\ref{thmc1}}
{\em  
Any two nearly cyclic 
pure isopairs 
are nearly unitarily equivalent
if and only if they have the same minimal polynomial.
}

In Section~\ref{sece}, we find a function-theoretic consequence of the operator theory.
Given a polynomial $q$, one can ask when $Y = Z_q \cap \T^2$ is polynomially convex. Apart from
the trivial case of when $q$ has factors of $(z - e^{i\theta})$
or $(w - e^{i\theta})$, the answer is that $Y$ fails to be polynomially convex if and only if
$q$ has an inner toral factor.

{\bf Theorem~\ref{thme1}}
{\em
Let $q$ be a polynomial in two variables with no linear factors. Then $Y = Z_q \cap \T^2$ is polynomially 
convex if and only if $q$ has no inner toral factor.
}

\section{Inner isopairs}
\label{secb}

\begin{definition} An \emph{isopair} is a pair $V=(V_1,V_2)$ of
  commuting isometries. An \emph{algebraic isopair} is an isopair that
satisfies a polynomial $p\in \mathbb{C}[z,w]$:
\[
p(V) = 0
\]
in which case $V$ may be called a $p$-isopair.
\end{definition}

\begin{definition} An isopair $V$ is \emph{pure} if 
\[
\bigcap_{m \geq 0} V_1^m  \mathcal{H} \= \{0\}
\=
\bigcap_{n\geq 0} V_2^n  \mathcal{H} .
\]
\end{definition}

Suppose $V = (V_1,V_2)$ is an isopair with $V_1$ pure. Let $k$ be the
dimension of the cokernel of $V_1$ (which is the Fredholm index of
$V_1^*$, and which we will call the multiplicity of $V_1$). Then
standard model theory for isometries, as described for example in
\cite{hel64} or \cite{szn-foi}, says that $V$ can be modeled on the
Hilbert space $H^2 \otimes \L$, where $\L$ is a Hilbert space of
dimension $k$, and $H^2$ is the Hardy space.  There is a $B(\L)$
valued inner function $\Phi$ so that $V$ is unitarily equivalent to
the pair $(M_z, M_\Phi)$, where $M_z$ is multiplication by the
independent variable (times $I_\L$) and $M_\Phi$ is multiplication by
the operator-valued function $\Phi$.  If $V_1$ is of finite
multiplicity, $k$ is finite and $\Phi$ is matrix-valued. If, in
addition, $V_2$ is of finite multiplicity, then $\Phi$ is a
matrix-valued rational inner function, \ie an analytic matrix-valued
function, each of whose entries is rational with poles outside the
closed unit disk, and such that everywhere on the unit circle the
matrix is unitary. 
Finally, if $V_2$ is also pure, this means that $\Phi$ is not the
direct sum of a constant unitary and another inner function. We shall
say in this case that the function $\Phi$ is pure.

\exam
A simple example is the pair $V = (V_1,V_2) = (M_{z^2},M_{z^3})$ on the classical
Hardy space $H^2$.  In this case $V$ satisfies $V_1^3-V_2^2=0$.  This
pair is unitarily equivalent to the pair $(M_z,M_\Phi)$ on $H^2\otimes
\C^2$ where
\[
\Phi(z) = \begin{pmatrix} 0 & z^2 \\ z & 0 \end{pmatrix}.
\]
The unitary equivalence $U$ comes from mapping $f = f_1(z^2)+zf_2(z^2) \in
H^2$ to $(f_1(z), f_2(z))^t \in H^2\otimes \C^2$, where here we are
dividing $f$ into its even and odd parts.  It is easy to check that
this Hilbert space isomorphism intertwines the two operator pairs:
\[
U M_{z^2} (f_1(z^2)+zf_2(z^2)) = U (z^2f_1(z^2)+z^3f_2(z^2)) =
\begin{pmatrix} zf_1(z) \\ zf_2(z) \end{pmatrix} = M_z \begin{pmatrix}
  f_1(z) \\ f_2(z) \end{pmatrix}
\]
(and similarly for $M_{z^3}$ on $H^2$ and $M_\Phi$ on $H^2\otimes
\C^2$).
\vs


Pure isopairs cannot be annihilated by arbitrary polynomials.
We shall show below (Theorem~\ref{thmb4}) that there is a minimal 
annihilating 
polynomial for any algebraic isopair. This 
minimal polynomial must be {\em inner toral.}
To define this, let us first establish the notation that $\D$ is the open unit disk, 
$\T$ is the unit circle, and $\E$ is the exterior of the closed unit disk in the plane.

\bd
A {\em distinguished variety} is an algebraic set $A$ in $\C^2$ such that 
$$
A \ \subseteq \ \D^2 \cup \T^2 \cup \E^2 .
$$
A polynomial is called {\em inner toral} if its zero set is a
distinguished variety.  \ed The terminology 
``inner toral'' 
(and explanation for it) 
is from
\cite{ams08}. The idea behind the name ``distinguished variety'' is
that the variety exits the bidisk through the distinguished boundary.
There is a close connection between pure algebraic isopairs and
distinguished varieties.  One theorem along these lines, proved first
in \cite{agmc_dv} and then, by a different method, in \cite{kn08ua},
is: \bt
\label{thmb1}
Let $A$ be a distinguished variety. 
Then there is a pure rational matrix-valued inner function $\Phi$ so that,
if
\be
\label{eqb01}
\det(\Phi(z) - wI) \= \frac{q(z,w)}{p(z,w)}, \ee 
then $A$ is the
zero-set of $q$.  Moreover, if $\Phi$ is any pure rational
matrix-valued inner function, and the polynomial $q$ is defined by
(\ref{eqb01}), then the zero set of $q$ is a distinguished variety.
\et
A sort of converse to Theorem~\ref{thmb1} is that the minimal annihilating polynomial
of any pure isopair is inner toral.

\begin{definition} $V = (V_1,V_2)$ is an \emph{inner isopair} if $V$
  is a pure isopair satisfying
\[
q(V) = 0
\]
where $q \in \mathbb{C}[z,w]$ is inner toral.
\end{definition}

\bt
\label{thmb2}
Every pure algebraic isopair is inner.
\et

To prove Theorem~\ref{thmb2}, we shall need some preliminary results.
First, we shall establish some basic facts about \emph{cyclic}
isopairs.

\begin{definition} An isopair $V$ is \emph{cyclic} if there exists $f \in
  \mathcal{H}$ such that
\[
\C[V] f:= \{p(V) f\ :\ p\inn \mathbb{C}[z,w]\}
\]
is dense in $\mathcal{H}$.
\end{definition}

\begin{definition}
We say a polynomial $p \in \C[z,w]$ has degree $(n,m)$ if it has
degree $n$ in $z$ and $m$ in $w$.
\end{definition}

\begin{lemma} \label{kernelprop} If $V$ is a cyclic isopair satisfying
  $p(V)=0$ where $p$ has degree $(n,m)$, then for all $(\alpha,\beta)
  \inn \D^2$, 
\beq
  \dim [\ker (V_1 - \alpha I)^* \cap \ker (V_2 - \beta I)^*] & \leq 1 \\
  \dim \ker (V_1 - \alpha I)^* & \leq m \\
  \dim \ker (V_2 - \beta I)^* & \leq n .  
\eeq
\end{lemma}
\bp
Applying appropriate M\"obius transformations to $V_1$ and $V_2$, we can assume
without loss of generality that $(\alpha, \beta) = (0,0)$.
Let $f$ be a cyclic vector.  

Suppose $g \in \ker V_1^* \cap \ker V_2^*$. Then
\[
\ip{Q(V) f}{g} = \ip{Q(0,0)f}{g}
\]
for any $Q\in \mathbb{C}[z,w]$.  If $\dim (\ker V_1^* \cap \ker V_2^*)
> 1$ then we could find a nonzero vector in $\ker V_1^* \cap \ker
V_2^*$ perpendicular to $f$. This would contradict cyclicity. So,
$\dim \ker V_1^* \cap \ker V_2^* \leq 1$.

If $\dim \ker V_1^* > m$ then we can choose $g \in \ker V_1^*$
perpendicular to $V_2^j f$ for $j=0,1,\dots, m-1$.  Now observe that 
\[
0 = p(V)^* g = p(0,V_2)^* g.
\]
Let $Q\in \mathbb{C}[z,w]$ and write
\[
Q(0,w) = s(w)p(0,w)+r(w) 
\]
where $r$ has degree less than $m$ (by the Euclidean algorithm).
Then, 
\beq
Q(V)^* g &= Q(0,V_2)^*g \\
&= s(V_2)^* p(0,V_2)^*g + r(V_2)^*g \\
&= r(V_2)^*g 
\eeq
and so
\[
\ip{Q(V) f}{g} = \ip{r(V_2)f}{g} = 0
\]
since $g$ is perpendicular to $V_2^j f$ for $j=0,1,\dots, m-1$.  This
contradicts cyclicity ($Q$ was arbitrary). So, $\dim \ker V_1^* \leq
m$.

Similarly, $\dim \ker V_2^*  \leq n$.
\ep

%


Let $\mult$ denote the multiplicity of an isometry:
\[
\mult V_i = \dim \ker V_i^*.
\]
We shall say that an isopair $V = (V_1,V_2)$ has finite multiplicity if both
$V_1$ and $V_2$ do.

\begin{lemma} 
\label{lemb6}
Let $V=(V_1,V_2)$ be a cyclic pure algebraic isopair and suppose $V$
is annihilated by an \emph{irreducible} polynomial $q \in
\mathbb{C}[z,w]$.  Then

 (1) $q$ is inner toral,

 (2) 
\[
\deg q = (\mult V_2, \mult V_1), \rm{\ and}
\]

(3) $q$ divides any polynomial $p$ that satisfies $p(V)
= 0$.
\end{lemma}

\bp  By Lemma~\ref{kernelprop}, 
we have
\be
\label{eqb8}
\deg q \geq (\mult V_2, \mult V_1) \ee (in each component separately).
Now, $V$ has a model as a pair of multiplication operators $(M_z,
M_\Phi)$ on $H^2\otimes \mathbb{C}^k$ where $k =\mult V_1$.  Since
$V_2$ has finite multiplicity (by Lemma~\ref{kernelprop}), $\Phi$ must
be a rational matrix valued inner function.

Let
\[
f(z,w_1,w_2) \= \frac{q(z,w_1) - q(z,w_2)}{w_1-w_2}.
\]
Letting
\be
\label{eqb23b}
Q(z,w) \= f(zI,wI,\Phi(z))   \= 
\left( q(zI,wI) - q(zI, \Phi(z)) \right) (wI - \Phi(z))^{-1},
\ee
we see that 
\be
\label{eqb9}
Q(z,w)(wI-\Phi(z)) = q(z,w)I.
\ee
Now, $Q$ is not
identically zero (else $q$ would be also), and $Q$ has lower degree in $w$ than
$q$.  So, the nonzero entries of the matrix polynomial $Q$ cannot
vanish identically on $Z_q$, the zero-set of $q$.  


As 
\beq
\left( wI - \Phi(z) \right) Q(z,w) &\=& q(z,w) I \\
&=& 0 \qquad {\rm\ on\ } Z_q ,
\eeq
we have 
\[
w \, Q(z,w) \= \Phi(z)\, Q(z,w) 
 \qquad {\rm\ on\ } Z_q .
\]

So
if $p \in \mathbb{C}[z,w]$ annihilates $V$, then 
\be
\label{eqb10}
p(zI,\Phi(z)) Q(z,w) \= p(z,w) Q(z,w) \= 0
 \qquad {\rm\ on\ } Z_q .
\ee
As $q$ is irreducible, and $Q$ does not vanish identically on $Z_q$, 
(\ref{eqb10}) shows that $q$ divides $p$, proving (3).

Now consider $d \in \C(z)[w]$, given by 
\[
d(z,w) = \det(wI-\Phi(z)).
\]
By Cayley-Hamilton, $d(z,\Phi(z))\equiv 0$, and therefore the numerator
of $d$ annihilates $V = (M_z,M_\Phi)$.  By the above, $q$ divides the
numerator of $d$ and since the degree of the numerator of $d$ is
$k=\mult V_1$, we see that $k$ is greater than or equal to the degree
of $q$ in $w$.  The reverse inequality is in (\ref{eqb8}), so these
two numbers are equal.  Interchanging the roles of $V_1$ and $V_2$, we
may conclude
\[
\deg q = (\mult V_2, \mult V_1)
\]
and this proves the second claim of the proposition.  Also, the zero
set of $d$ is inner toral and this implies $Z_q$ is inner toral, since
$Z_q \subset Z_d$.  This proves the first claim.
\ep

If $V$ is a cyclic isopair, then in particular it is a cyclic
subnormal pair, and so has another nice representation. For $\mu$ any
compactly supported measure in $\C^2$, let $\ptm$ denote the closure
of the polynomials in $\ltm$. Then we have the following
representation; see \cite{con88, mcc93a} and references therein for
details.

\bt
Let $V$ be a cyclic isopair on the Hilbert space $\h$, with cyclic vector $u$. Then there is a
positive Borel measure $\mu$ on $\T^2$ and a unitary operator $U$ from $\h$ onto $\ptm$ that
maps $u$ to the constant function $1$, and such that $U$ intertwines $V$ with the pair $(M_z,M_w)$
of multiplication by the coordinate functions.
\label{thmb3}
\et

Theorem~\ref{thmb3} makes it easy to prove that the minimal polynomial
of a pure algebraic isopair is square-free.

\begin{lemma}
\label{lemb9}
Suppose $V$ is a pure $p$-isopair, and the irreducible factors of $p$ are $p_i$, each with multiplicity $t_i$:
\[
p \= \prod p_i^{t_i} .
\]
Let $q = \prod p_i$. Then $q(V) = 0$.
\end{lemma}
\bp
Choose some vector $u$. Let 
\[
\mathcal{K} \= \overline{\mathbb{C}[V]u} 
\]
and
let $T = V\mid_{\mathcal{K}}$.
By Theorem~\ref{thmb3}, $T$ is unitarily equivalent to $(M_z,M_w)$ on some $\ptm$.
As $p(T) = 0$, we must have that $p$ vanishes on the support of $\mu$.
Therefore so does $q$, and so
\[
\| q(V) u \|^2 \= \int |q|^2 d\mu \= 0.
\]
As $u$ was arbitrary, we must have that $q(V) = 0$.
\ep

\begin{lemma} 
\label{lemb7}
Suppose $V$ is a pure $q$-isopair where $q \in
  \mathbb{C}[z,w]$ is a product of distinct irreducible factors and
  $V$ is not annihilated by any factor of $q$.  Then, $q$ is inner toral
  and divides any polynomial that annihilates $V$.
\end{lemma}

\bp First, we claim that any irreducible factor of $q$ is
  inner toral.  Let $q_0$ be an irreducible factor and write $q = q_0
  q_1$.  Then, $u:=q_1(V) u_0 \ne 0$ for some $u_0$ in our Hilbert
  space.  Consider the cyclic subspace $\mathcal{K}$ generated by $u$
\[
\mathcal{K} \ :=\  \overline{\mathbb{C}[V]u} \= \vee \{g(V)u: g \in
\mathbb{C}[z,w]\},
\]
where $\vee$ denotes the closed linear span. 
Let $T = V\mid_{\mathcal{K}}$ be the pure $q_0$-isopair obtained by
restricting $V$ to the invariant subspace $\mathcal{K}$.  By Lemma~\ref{lemb6},
$q_0$ must be inner toral.  As $q_0$ was an arbitrary
irreducible factor, all factors of $q$ are inner toral, and this
implies $q$ is inner toral.

Also, if $g(V)=0$ for some $g\in \mathbb{C}[z,w]$, then $g(T) = 0$ and
by Lemma \ref{lemb6} this implies $q_0$ divides $g$. As $q_0$ was an
arbitrary irreducible factor of $q$, we see that $q$ divides $g$.
This proves the second claim of the lemma. 
\ep

Putting together what we have proved, we get the following theorem,
which contains Theorem~\ref{thmb2}.

\bt
\label{thmb4}
Let $V$ be a pure algebraic isopair. Then there exists a square-free
inner toral polynomial $q$ that annihilates $V$.  Moreover, if $p$ is
any polynomial that annihilates $V$, then $q$ divides $p$.  
\et

We shall call the polynomial $q$ the {\em minimal polynomial} of $V$.

\section{Decomposition of algebraic isopairs} \label{secba0}

In this section we show how algebraic isopairs are {\em nearly} a
direct sum of algebraic isopairs annihilated by irreducible polynomials. 

\bt
\label{thmpropb5} Let $V$ be a pure algebraic isopair with
  minimal polynomial $q$, and let $q_1,q_2, \dots, q_N$ be the
  (distinct) irreducible factors of $q$.  If $V$ has finite
  multiplicity, then $V$ has a finite codimension invariant subspace
  $\mathcal{K}$ on which
\[
V\mid_{\mathcal{K}} = W_1 \oplus W_2 \oplus \cdots \oplus W_N
\]
where $W_j$ is a $q_j$-isopair, $j=1,\dots, N$.
\et

As a quick example, consider the reducible algebraic set $z^2=w^2$. We can
define a Hilbert space by defining
\[
||p||^2 = \int_{0}^{2\pi} |p(e^{i\theta}, e^{i\theta})|^2 d\theta +
\int_{0}^{2\pi} |p(e^{i\theta}, -e^{i\theta})|^2 d\theta
\] 
for each $p \in \C[z,w]$ and then completing this to an $H^2$ space.
The pair $(M_z,M_w)$ will be a $(z^2-w^2)$-isopair and the decomposition
from the above proposition consists of letting $\mathcal{K}$ be the 
functions that vanish at $(0,0)$, and dividing this Hilbert space into
the functions that are a multiple of $z-w$ and  those that are a multiple
of $z+w$. $\K^\perp$ is the constant functions.

Before we prove the theorem we need the following.

\begin{lemma} \label{lemba01} Suppose $p \in
  \C[z,w]$ is inner toral and reducible $p=p_1p_2$.  If $V$ is a pure
  $p$-isopair, then
\[
{\rm ran\ }p_1(V) \perp {\rm ran\ }p_2(V)
\]
(ran denotes the range).
\end{lemma}

\bp
Since $p_1$ is inner toral, it is a fact that $p_1$ is symmetric in
the sense that
\[
z^nw^m\overline{p_1(1/\bar{z},1/\bar{w})} = \mu p_1(z,w)
\]
where $\mu$ is a unimodular constant and $(n,m)$ is the degree of $p_1$
(see \cite{kn08ua}). In fact, we may assume $\mu=1$ by replacing $p_1$
with an appropriate constant multiple.  Hence, if we write
\[
p_1(z,w) = \sum_{j=0}^{n}\sum_{k=0}^{m} a_{jk} z^j w^k
\]
it follows that $\overline{a_{jk}} = a_{(n-j)(m-k)}$.  This can be
used to deduce
\[
p_1(V)^* V_1^n V_2^m = p_1(V)
\]
since $V_1,V_2$ are isometries.

Then,
\beq
p_1(V)^* p_2(V) &= p_1(V)^* (V_1^n V_2^m)^* V_1^n V_2^m p_2(V) \\
&= (V_1^n V_2^m)^* p_1(V)^* V_1^n V_2^m p_2(V) \\
&= (V_1^n V_2^m)^* p_1(V) p_2(V) \\
&= (V_1^n V_2^m)^* p(V) = 0
\eeq
So, for any $f,g \in \mathcal{H}$
\[
\ip{p_2(V)f}{p_1(V) g} = \ip{p_1(V)^*p_2(V) f}{g} = 0.
\]
Hence,
\[
\ran p_1(V) \perp \ran p_2(V)
\]
\ep

{\sc Proof of Theorem~\ref{thmpropb5}:}
  Let $p = q_2q_3\cdots q_N$ (where $q_1,q_2, \dots, q_N$ come from
  the statement of the proposition).  We will show that $V$ has a finite
  codimension invariant subspace on which $V$ is the direct sum of a
  $q_1$-isopair and a $p$-isopair.  The proposition will then follow
  by induction.

  By Lemma \ref{lemba01}, $\ran q_1(V)$ and $\ran p(V)$ are
  orthogonal.  Let $\mathcal{K}' = (\ran q_1(V)+\ran p(V))^{\perp}$.
  The assumption that $V$ has finite multiplicity implies that
  $\mathcal{K}'$ is finite dimensional as follows. First note that as
  in the previous lemma we may assume that $p$ and $q_1$ are
  symmetric, so that we have the formulas:
\[
p(V)^* = V_1^{*n}V_2^{*m}p(V) \qquad q_1(V)^* = V_1^{*j}V_2^{*k}q_1(V)
\]
where the degree of $p$ is $(n,m)$ and the degree of $q_1$ is $(j,k)$.
If $f \in \mathcal{K}'$, then $0=p(V)^*f = V_1^{*n}V_2^{*m}p(V)f$ and
$0=q_1(V)^*f= V_1^{*j}V_2^{*k}q_1(V)$.  Since $V$ is assumed to have
finite multiplicity, the kernels of $V_1^{*n}V_2^{*m}$ and
$V_1^{*j}V_2^{*k}$ are both finite dimensional.  Hence, the ranges of
$p(V)\mid_{\mathcal{K}'}$ and $q_1(V)\mid_{\mathcal{K}'}$ are both finite
  dimensional (since they map into ${\rm ker} V_1^{*n}V_2^{*m}$ and
  ${\rm ker}V_1^{*j}V_2^{*k}$ respectively).

Now, since $q_1$ and $p$ are relatively prime, there exist nonzero polynomials
$A, B \in \C[z,w]$, $C \in \C[z]$ such that
\[
A(z,w)q_1(z,w)+B(z,w)p(z,w) = C(z)
\]
(let $C$ be the resultant of $q_1$ and $p$).  Substituting $V$
\[
A(V)q_1(V)+B(V)p(V) = C(V_1)
\]
it is then apparent that $C(V_1)\mid_{\mathcal{K}'}$ has finite
dimensional range.  If $\mathcal{K}'$ were infinite dimensional then
$C(V_1)$ would have nontrivial kernel.  This is impossible (a pure
isometry cannot have eigenvalues for instance), so $\K'$ is finite
dimensional.  

It is clear that $\mathcal{K}_1 = \ran p(V)$ and $\mathcal{K}_2 =
\ran q_1(V)$ are mutually orthogonal invariant subspaces for $V$.
Also, $V\mid_{\mathcal{K}_1}$ is a $q_1$-isopair and
$V\mid_{\mathcal{K}_2}$ is a $p$-isopair.  Since $\mathcal{K}_1\oplus
\mathcal{K}_2$ has finite codimension, the proposition is proved with
$\mathcal{K} = \mathcal{K}_1\oplus \mathcal{K}_2$.
\ep

\section{Nearly cyclic isopairs}
\label{secc}

\bd
An isopair $V$ is {\em nearly cyclic}
if there is a vector $u$ such that 
\[
\overline{\mathbb{C}[V]u} \= \vee \{g(V)u: g \in
\mathbb{C}[z,w]\}
\]
is of finite codimension.
\ed

For example, the pair $(M_{z^2}, M_{z^3})$ on $H^2(\T)$, is not
cyclic, because for any $f \in H^2(\T)$, we can find a function $g$
orthogonal to $\C[z^2,z^3]f$.  Namely, write 
\[
f(z) = a + bz + \rm{\ higher\ order\ terms}
\]
and define $g(z) = -\bar{b}+\bar{a}z$.  Then, $\ip{f}{g}=0$ and since
$g$ is linear it is orthogonal to multiples of $z^2$.  On the other
hand, the pair $(M_{z^2},M_{z^3})$ is {\em nearly} cyclic:
\[
\overline{\mathbb{C}[M_{z^2},M_{z^3}]1} = H^2\ominus \C \{z\}.
\]

\bd
Two isopairs $V = (V_1,V_2)$ on $\h$ and $V' = (V_1',V_2')$ 
on $\h'$ are {\em nearly unitarily equivalent} if
there is a finite codimension 
$V$-invariant subspace $\K$ of $\h$, a 
finite codimensional 
$V'$-invariant subspace $\K'$ of $\h'$, and unitary
operators $U: \h \to \K'$ and $U' : \h' \to \K'$ such that
\beq
 U V_r U^* &\=& V_r' |_{\K'} \qquad r =1,2 \\
 U' {V'}_r {U'}^* &\=& V_r |_{\K} \qquad r =1,2 .
\eeq
(In words, each one is unitarily equivalent to the other restricted to a finite codimensional
invariant subspace).
\ed 
The principal result of this section is:
\bt
\label{thmc1}
Any two nearly cyclic 
pure isopairs 
are nearly unitarily equivalent
if and only if they have the same minimal polynomial.
\et

For the rest of this section, fix some square-free inner toral polynomial $q$.

The necessity is obvious, as restricting an algebraic isopair to a 
finite codimensional
invariant subspace
does not change the minimal polynomial.
Let us give an overhead view of the proof of sufficiency in Theorem~\ref{thmc1}.  To
show any two nearly cyclic pure algebraic isopairs with the same
annihilating polynomial are nearly equivalent it suffices to show (1)
that any such nearly cyclic isopair has a finite dimensional extension
to a cyclic isopair (and therefore the nearly cyclic isopairs can be
extended \emph{and} restricted to cyclic ones) and (2) that all such
cyclic isopairs are nearly unitarily equivalent to one particular
choice of a nearly cyclic isopair $W$.  The idea is if $V$ is a given
nearly cyclic isopair then we have the following diagram
\[
W \to {\rm\ some\ cyclic\ } \to V \to {\rm some\ cyclic\ } \to W
\]
where the arrows denote restrictions to finite codimensional invariant
subspaces.  

The study of cyclic isopairs and the construction of $W$ requires us
to lift many of our questions to a finite Riemann surface that
desingularizes $Z_q\cap \D^2$.  Let $\O = Z_q \cap \D^2$.
As described in \cite{amhyp}, there is a finite Riemann surface $S$
and a holomap $h$ from $S$ onto $\O$ (a holomap is a proper
holomorphic map that is one-to-one and non-singular except on finitely
many points).  We shall let $A(S)$ denote the algebra of functions
that are holomorphic on $S$ and continuous up to the boundary, and we
shall let $A_h(S)$ be the finite codimensional subalgebra that is the
closure of polynomials in $h = (h_1,h_2)$.  Note that $A_h(S)$ can be
described by a finite number of (homogeneous) linear relations on
derivatives of elements of $A(S)$ at a finite number of specified
points.  This implies that we can find an element of $g \in A_h(S)$
that multiplies $A(S)$ into $A_h(S)$:
\[
g A(S) \subset A_h(S).
\]
We simply must choose $g$ to vanish to sufficient order at a finite
number of points.

Let $\om$ be harmonic measure on $S$ at some point
$z_0$, fixed hereinafter.  Let $\ato$ be the closure in $L^2(\om)$ of $A(S)$ and let $W$
be $M_h=(M_{h_1},M_{h_2})$ on $\ato$.  Our eventual goal is to show
that every cyclic isopair is nearly unitarily equivalent to $W$.

In addition, we can elaborate on the structure of cyclic isopairs via
$h$.  Let $V$ be a cyclic $q$-isopair.  We may model $V$ as
multiplication by coordinate functions $(M_z,M_w)$ on $\ptm$ for some
measure $\mu$ on $Z_q\cap \T^2$.  Let $\nu = h^*(\mu)$ be the
pull-back of $\mu$ to $X:= \partial S$ (\ie $\nu(E) = \mu(h(E))$ for
any Borel subset of $X$).  Let $\atn$ be the closure in $L^2(\nu)$ of
$A(S)$, and let $\ath$ be the closure of $A_h(S)$ in $L^2(\nu)$.

Then $V$ is unitarily equivalent to $M_h = (M_{h_1}, M_{h_2})$ on
$\ath$.  As $\ath$ is of finite codimension in $\atn$, $V$ has a
finite-dimensional extension to an isopair $V^S$ that is unitarily
equivalent to $M_h$ on $\atn$. Let us record these observations and a
few more.

\bl 
\label{lemcd}
If $V$ is a cyclic pure $q$-isopair, then there exists a positive
Borel measure $\nu$ on $X$ such that $V$ is unitarily equivalent to
$M_h$ on $\ath$.  Furthermore, the measures $\nu$ and $\om$ are
mutually absolutely continuous and satisfy
\[
\int \log \frac{d\nu}{d\om} d\om > -\infty.
\]
\el \bp By a result of J.~Wermer \cite{wer64}, the algebra $A(S)$ is a
hypo-Dirichlet algebra on $X=\partial S$ , \ie the real parts of
functions in $A(S)$ form a finite codimensional subspace of $C_\R(X)$.
As $V$ is pure, $\mu$ is non-atomic (an atom would yield a bounded
point evaluation on $Z_q\cap\T^2$ and the corresponding evaluation
kernel would be an eigenvector for $V^*$). Hence $\nu$ is non-atomic
and $V^S$ is also pure. (Note that
a finite dimensional extension of a pure
pair can only fail to be pure if the extension has a unitary summand
and hence an eigenvalue.  However, in this case we can multiply
elements of $\atn$ by some $g\in A_h(S)$ and produce an element of
$\ath$.  An eigenvector $f$ of $V^S$ would then produce an eigenvector
$gf$ of $V$, contradicting purity of $V$.)

Therefore $\atn \neq L^2(\nu)$,
and $\atn$ has no $L^2$ summand. 
Next, we claim that $\nu$ is absolutely continuous with respect to $\om$,
using an argument from \cite{mcc90b}.
Indeed, let us write $\nu_a$ and $\nu_s$ for the absolutely continuous and
singular parts of $\nu$. Let $E$ be an $F_\sigma$ 
set such that $\nu_a(E) = 0$ and $\nu_s( X \setminus E) = 0$.
By Forelli's lemma (\cite[II.7.3]{gam}, applicable here because of
Corollary 1 to Theorem 3.1 of \cite{ahesar1}), there is a sequence
$f_n$ in $A(S)$ with $\| f_n \|_X \leq 1$, and such that
$f_n(x)$ tends to $0$ for every $x$ in $E$, and to $1$ 
for $\om$-a.e. $x$.
Some subsequence of $f_n$ converges weak-* to a function $g$.
By the dominated convergence theorem, $\int g d\om = 1$, so
$g = 1\ \om$-a.e. Again using dominated convergence, we see that
for every $h$ in $L^2(\nu_s)$ we have $\int g h d \nu = 0$.
Therefore $1-g$, which is in $A^2(\nu) \cap L^\i(\nu)$, agrees
with the characterictic function of $E$ $\nu$-a.e., so
$$
A^2(\nu) \= A^2(\nu_a) \oplus A^2(\nu_s) .
$$
By the Kolmogorov-Krein theorem \cite[V.8.1]{gam}, 
$A^2(\nu_s) = L^2(\nu_s)$, so we conclude that $\nu_s$ is null,
as claimed.

As $\atn \neq L^2(\nu)$, it follows from 
the work of P.~Ahern and D.~Sarason on hypo-Dirichlet algebras
\cite[Corollary to Thm. 10.1]{ahesar1}, that
\be
\label{eqc2}
\int \log \frac{d\nu}{d\om} \, d\om  \ > \ -\i.
\ee

%
%

\ep

The next proposition allows us to dispense with dealing with
\emph{nearly} cyclic isopairs.

\begin{prop} Any nearly cyclic pure $q$-isopair is unitarily equivalent to
  a cyclic pure $q$-isopair restricted to a finite codimensional
  invariant subspace.
\end{prop}

\bp Let $V$ be a nearly cyclic pure $q$-isopair on the Hilbert space
$\mathcal{H}$, let $\mathcal{K}$ be a finite codimensional invariant
subspace on which $V$ is cyclic, and let $\mathcal{F} = \mathcal{H}
\ominus \mathcal{K}$ (a finite dimensional subspace).  Since
$V\mid_{\mathcal{K}}$ may be modeled as $(M_z,M_w)$ on $\ptm$ for some
measure $\mu$ supported on $Z_q\cap \T^2$, we shall simply identify
$\mathcal{K} = \ptm$ and $V\mid_{\mathcal{K}} = (M_z,M_w)$.  Then, the
pair $V$ can be written in block form as
\[
V = \begin{array}{cc} & \begin{array}{cc} \ptm &
    \mathcal{F} \end{array} \\
\begin{array}{c} \ptm \\ \mathcal{F} \end{array} & \begin{pmatrix}
  (M_z,M_w) & (B_1,B_2) \\ 0 & (A_1,A_2) \end{pmatrix} \end{array}
\]
where $(A_1,A_2)$ is a pair of commuting contractions on the finite
dimensional space $\mathcal{F}$ that (by purity) have no unimodular
eigenvalues.

Let $u:\D\to\D$ be a finite Blaschke product that annihilates $A_1$:
\[
u(A_1) = 0.
\]
If we apply such a $u$ to $V_1$ (we can do this since $u$'s power series
is absolutely convergent in $\overline{\D}$) we get an isometry
\[
u(V_1) = \begin{array}{cc} & \begin{array}{lr} \ptm &
    \mathcal{F} \end{array} \\
\begin{array}{c} \ptm \\ \mathcal{F} \end{array} & \begin{pmatrix}
  M_u & u(B_1) \\ 0 & 0 \end{pmatrix} \end{array}
\]
where $M_u$ is multiplication by $u$ on $\ptm$.  In particular, the
range, say $\mathcal{L}$, of $u(V_1)$ is contained in $\ptm$.
Therefore, $u(V_1): \mathcal{H} \to \mathcal{L}$ can be thought of as a
Hilbert space isomorphism that intertwines $V$ on $\mathcal{H}$ and
$(M_z,M_w)$ on $\mathcal{L} \subset \ptm$.  This proves $V$ can be
modeled as a restriction of $(M_z,M_w)$ on $\ptm$ (\ie a cyclic pair)
to an invariant subspace (\ie $\mathcal{L}$).  The key thing left to
prove is that $\mathcal{L}$ has finite codimension in $\ptm$.  For
this it suffices to prove $u\ptm$ has finite codimension in $\ptm$
since $u\ptm \subset \mathcal{L}$.

To see this, it helps to use the $\atn$ model described above (in
Lemma \ref{lemcd}).  Namely, we need to prove $u(h_1)\ath$ has finite
codimension in $\ath$.  Since $\ath$ has finite codimension in $\atn$,
it suffices to prove $u(h_1)\atn$ has finite codimension in $\atn$.
This follows from the fact that $u(z)$ vanishes finitely often among
$(z,w) \in Z_q$ (and therefore $u(h_1)$ has finitely many zeros on
$S$) and so any analytic function $f$ in $A(S)$ which vanishes to
higher order at $u(h_1)$'s zeros than $u(h_1)$, is a multiple of
$u(h_1)$: $f/(u(h_1)) \in A(S)$.  \ep

Thus, any nearly cyclic has an finite dimensional extension to a
cyclic and by definition a finite codimension restriction to a cyclic.
So, now we let $V$ be a pure cyclic $q$-isopair, which we think of as
$M_h$ on $\ath$ and we let $V^S$ be the extension of $V$ to
$\atn$. The following lemma is really a chain of lemmas, since we
prefer to introduce ideas from references as we need them. Recall that
$W$ refers to $M_h$ on $\ato$.

\bl
\label{lemc1}
With notation as above, $V$ is unitarily equivalent to $W$ restricted
to a finite codimensional invariant subspace.
\el

\bp
Suppose we can find 
a function $f$ in $\ato$ that is {\em nearly outer},
in the sense that the invariant subspace it generates,
$$
[f] \ := \
\overline{A(S) \,f} ,
$$
 is of finite 
codimension, and such that $\dis |f|^2 = \frac{d\nu}{d\om} $.
Let
$$
[f]_W \ := \
\overline{\mathbb{C}[W]f} \= 
\overline{A_h(S) \,f} ,
$$
which will be of finite codimension in $[f]$.

Then the map $ g\cdot f \mapsto g$ extends to a unitary between
$[f]_W$ in $\ato$ and $\athn$ that intertwines $W |_{[f]}$ and $V$.
By Lemma~\ref{lemc2}, which we prove below, such a nearly outer
function exists.  \ep

When does a nearly outer function exist with a given log-integrable
modulus?  Let $L$ be the codimension of $Re(A(S))$ in $C_\R(X)$.

Ahern and Sarason proved that any log-integrable positive function can
be written as $|f|^2$ for some $f$ in $\ato$, and they conjectured that
$f$ can be chosen so that $[f]$ is of codimension no more than the
codimension of $\Re(A(S))$ in $C_\R(X)$ \cite{ahesar1}. This
conjecture is still open, though it has been proved in the planar case
by G.~Tumarkin and S.Ya.~Khavinson \cite{tumkha58}.

However, using results of S.Ya.~Khavinson \cite{kh89a,kh89b} for general finite Riemann surfaces,
which generalize results of D.~Khavinson for the planar case \cite{kh84,kh85}, we can prove that
$f$ can be chosen with $[f]$ of finite codimension.

The idea is that when we write down the Green integral (or Poisson
integral) of a measure, the resulting harmonic function's conjugate
function will in general be multi-valued.
To fix this, we need to worry about the periods on a homology basis
$K_r$ for $S$.  There are $L = 2h + n - 1$ such curves, where $S$ has
$h$ handles and $n$ boundary components.  Choose $L$ disjoint arcs
$\Delta_j$ in $\partial S$, and positive measures $\nu_j$ supported on
each arc, so that each $\nu_j$ is boundedly absolutely continuous with
respect to $\om$ (in fact we may simply take $\nu_j$ to be harmonic
measure restricted to $\Delta_j$).  Khavinson shows that after
shrinking $\Delta_j$ if necessary, the matrix $A$ of periods of the
harmonic conjugate of the Green integrals of the $\nu_j$ along the
curves $K_r$ is non-singular \cite{kh89a}, and hence may be used to
``correct'' the periods of other functions (while at the same time we
have some control over what is happening on the boundary).

The Green kernel is defined by
$$
P(z,\zeta) 
\= \frac{1}{2\pi} 
\frac{ \frac{\partial}{\partial n_\zeta} G(z, \zeta)}
{\frac{\partial}{\partial n_\zeta} G(z_0, \zeta)}
\qquad z \inn S,\ \zeta \inn \partial S ,
$$
where $G(z,\zeta)$ is the Green's function 
with pole at $z$, and $n_\zeta$ is the outward normal.  The Green
integral of a measure $\nu$ is then
\[
\int_{\partial S} P(z,\zeta) d\nu(\zeta).
\]
(Note we are using $z,\zeta$ to refer to points of $S$, while these
letters are typically reserved for uniformizers on Riemann surfaces.)


Let 
$$
\omega_j(z) \= \int_{\Delta_j} P(z,\zeta) d\nu_j(\zeta) .
$$

Let us explicitly define:

\begin{definition} A function $f \in \ato$ is \emph{nearly outer} if
\[
[f] := \overline{A(S)f}
\]
has finite codimension in $\ato$.
\end{definition}

\bl
\label{lemc2}
Any log-integrable function $w$ on $\partial S$ is the modulus squared
of a nearly outer function.  
\el 
\bp
First, we write down the Green integral of $1/2\ \log w$
\[
\int_{\partial S} P(z, \zeta) \frac{1}{2} \log w(\zeta) d\om(\zeta).
\]
There exist real constants $\lambda_j$ such that
\begin{equation} \label{hfunc}
h(z) = \int_{\partial S} P(z, \zeta) \, [ \frac{1}{2} \log w(\zeta) d\om(\zeta) -
\sum_{j=1}^{L} \lambda_j d\nu_j(\zeta) ]
\end{equation}
has a single-valued harmonic conjugate $*h$, and $h$ has boundary
values given by
 \[
 h(\zeta) = \begin{cases} \frac{1}{2}\log w(\zeta) & \zeta
   \in \partial
   S \setminus \cup\Delta_j \\
   \frac{1}{2} \log w(\zeta) -\lambda_j \frac{d\nu_j}{d\om} (\zeta)&
   \zeta \in \Delta_j
 \end{cases}
 \] 
\ie
\begin{equation} \label{hfunc2}
h d\om = (1/2)\log w d\om - \sum_j \lambda_j d\nu_j.
\end{equation}

The function
\[
g(z) = \exp(h(z)+i*h(z))
\]
is outer in the sense that
\[
\log|g(z_0)| = \frac{1}{2\pi} \int_{\partial S} \log|g(\zeta)| d\om(\zeta)
\]
(by \eqref{hfunc} and \eqref{hfunc2} since $P(z_0,\zeta)=
\frac{1}{2\pi}$).


It follows from \cite{ahesar1} (Theorem 7.1 and the discussion
following Theorem 9.1) that it is also outer in the sense that $[g]$
is all of $\atn$.  If we can find a finite Blaschke product $F$ whose
modulus on the boundary is \be
\label{eq3}
\log |F(\zeta)| \=
\left\{ 
\begin{array}{cl}
0 & \zeta \inn \partial S \setminus \cup \Delta_j \\
\lambda_j \frac{d\nu_j}{d\om}(\zeta)& \zeta \inn \Delta_j ,
\end{array}
\right. , \ee then $f = Fg$ will be a nearly outer function ($g$ is
outer and any function that vanishes on the zeros of $F$ will be
divisible by $F$) satisfying $|f|^2 = w$ a.e. on $\partial S$.  We
shall prove $F$ exists in Lemma~\ref{lemck}.  \ep

Following Schiffer-Spencer \cite{ss54},
Khavinson defines a basis for the space of abelian differentials of
the first kind via the Green's function.  Specifically, $dZ_j$ is
defined using the local expression
\[
Z_j'(z) = -\frac{1}{\pi} \int_{K_j} \frac{\partial^2
  G(z,\zeta)}{\partial z \partial \zeta} d\zeta
\]
or 
\[
\Im Z_j(z) = 2i \frac{1}{2\pi} \int_{K_j}
\frac{\partial G(z,\zeta)}{\partial \zeta} d\zeta
\]
where again $K_1, \dots, K_L$ form a canonical homology basis for $S$.
It should be noted that ${\Im Z_j}$ is single valued on $S\setminus
K_j$ (and hence single valued everywhere when $K_j$ is a boundary
cycle), but has a jump across $K_j$ when $K_j$ is a cycle
corresponding to a handle. 

\begin{lemma} 
\label{lemck}
Given real constants $\lambda_j$, $j=1,\dots,L$, there exists
  a bounded holomorphic function $F:S\to \C$ with finitely many zeros
  in $S$ and boundary modulus satisfying
\be \label{eq3-2}
\log |F(\zeta)| = \begin{cases} 0 & \zeta \inn \partial S\setminus
  \cup\Delta_j \\
\lambda_j\frac{d\nu_j}{d\om}(\zeta) & \zeta \in \Delta_j \end{cases}
\ee
\end{lemma}

\bp
Khavinson shows that for each point $\alpha$ in $S$ there is a function $B(z;\alpha)$ 
that has a single zero at $\alpha$, and 
such that 
\be \label{Bvalues}
\log |B(z;\alpha)| \=
\left\{ 
\begin{array}{cl}
0 & z \inn \partial S \setminus \cup \Delta_j \\
c_j \frac{d\nu_j}{d\omega} & z \inn \Delta_j 
\end{array}
\right. .
\ee
Moreover, one finds the $c_j$ (real constants) by using the
(invertible) period matrix $A$ and the abelian differentials $dZ_j$ by
the formula
$$
\left(
\begin{array}{c}
c_1\\
\vdots\\
c_L
\end{array}
\right)
\=
- 2 \pi A^{-1}
\left(
\begin{array}{c}
\Im Z_1(\alpha)\\
\vdots\\
\Im Z_L(\alpha)
\end{array}
\right) .
$$

Also, for each $(d_1,\dots,d_L)^t \in \Z^L$, if we set
\[
\begin{pmatrix} c_1 \\ \vdots \\ c_L \end{pmatrix} = -2\pi
A^{-1} \begin{pmatrix} d_1 \\ \vdots \\ d_L \end{pmatrix}
\]
then there is a holomorphic function $B$ on $S$ with \emph{no zeros}
and boundary modulus values satisfying equation \eqref{Bvalues}.  

Therefore to prove the lemma it suffices to show that we can take a
finite positive integer combination of vectors of the form $(\Im
Z_1(\alpha), \dots, \Im Z_L(\alpha))^t$ and obtain every element of
$\R^L/\Z^L$, for then we could find a combination satisfying 
\[
-\frac{1}{2\pi} A \begin{pmatrix} \lambda_1 \\ \vdots
  \\ \lambda_L \end{pmatrix} = \sum_{j=1}^k \begin{pmatrix} \Im Z_1(\alpha_j) \\ \vdots \\ \Im
  Z_L(\alpha_j) \end{pmatrix} \qquad {\rm mod\ } \Z^L
\]
and upon modifying it by an element $(d_1,\dots,d_L)^t \in \Z^L$, we
would obtain a finite Blaschke product with zeros at the $\alpha_1,
\dots, \alpha_k$ and the desired boundary modulus values.  For this it
suffices to prove the following claim.

Claim: As $\alpha_1, \dots, \alpha_L$ vary over $S$,
the vectors
\be
\label{eq1}
\left(
\begin{array}{c}
\sum_{r=1}^L \Im Z_1(\alpha_r)\\
\vdots\\
\sum_{r=1}^L \Im Z_L(\alpha_r)\\
\end{array}
\right) 
\ee
have interior in $\R^L$.

Indeed, given the claim, it follows that there is some finite $N$ so
that sums of $N$ vectors of the form (\ref{eq1}) (\ie with $LN$ points
$\alpha_r$) form a ball large enough that
it covers an entire
cell of $\R^L / \Z^L$, and so a 
Blaschke product of degree $LN$ will satisfy
(\ref{eq3-2}).

Proof of claim: There is no harm in assuming our argument takes place
inside some coordinate neighborhood.  Consider the derivative of
(\ref{eq1}) with respect to $\alpha_1, \dots, \alpha_L$.  
If the
$L$-by-$L$ matrix $ ( Z_j'(\alpha_r) )$ is of full rank, 
then for some choice of unimodular $\tau_r$ the real matrix $ (\Im
(\tau_r Z_j'(\alpha_r)) )$ is invertible (a linear algebra exercise).  
This proves the claim in this case because if we replace $\alpha_r$ by
$\tau_r \alpha_r$ in \eqref{eq1} and take derivatives we get a
nonsingular Jacobian matrix.

Otherwise, there are real numbers $c_j$, not all zero, so that, for
every $1 \leq r \leq L$, we have \be
\label{eq54}
\sum_{j=1}^L c_j Z_j'(\alpha_r) = 0.  \ee But (\ref{eq54}) means that
the differential $\sum c_j dZ_j$, which extends by reflection to the
double of $S$, vanishes at $L$ points in $S$ and $L$ more on the
reflection.  The double has genus $L$, and so a differential of the
first kind must have exactly $2L-2$ zeroes (see \cite{ss54}
(3.5.1)). As the $dZ_j$'s are linearly independent, this forces all
the $c_j$'s to be zero, a contradiction.  \ep

This chain of lemmas completes the proof of Lemma \ref{lemc1}, which
says that any cyclic $V$ is a finite codimensional restriction of $W$.

\vs

{\sc Proof of Theorem~\ref{thmc1}}. In light of Lemma~\ref{lemc1}, it
remains to show that $V$, a cyclic $q$-isopair, can be restricted to a finite
codimensional subspace to become unitarily equivalent to $W$.  Again we view
$V$ as $M_h$ on $\ath$.  Now $A_h(S)$ has finite codimension in $A(S)$
and is defined by a finite number of linear relations on derivatives
(at a finite number of points) of elements of $A(S)$ (see
\cite{amhyp}).  In particular, any element of $A^2(\nu)$ that vanishes
to high enough order at these finite points will be inside $\ath$.
This can be accomplished by multiplying $\atn$ by an appropriate
finite Blaschke product $F$.  Then, $F\atn \subset \ath$ and the
operator $M_h$ restricted to $F\atn$ is unitarily equivalent to $M_h$
on $A^2(\sigma)$ where $\sigma = |F|^2\nu$.  This proves $V$ has a
finite codimension restriction to $M_h$ on $A^2(\sigma)$.


So it suffices to show that one can 
find a unitary equivalence between
$M_h$ on $\ats$ restricted to an
invariant subspace of finite codimension and
$M_h$ on $\ato$.
But this can be done by finding a nearly outer function $f$ with modulus
$|f|^2 = \frac{d\omega}{d\sigma}$ just as in
Lemma~\ref{lemc2}.
\ep
\vs
\rem
\label{remc7}
In the proof of Theorem~\ref{thmc1}, we did not strongly use the fact that $V$
and $W$ are isopairs. We could more generally look at nearly cyclic pure subnormal pairs whose
spectral measures were supported on the boundary of some hyperbolic algebraic set.

\vs

We can translate Theorem~\ref{thmc1} into the matrix models and get the
following result.

\begin{cor} Suppose $(M_z,M_\Phi)$ and $(M_z,M_\Psi)$ are two
  nearly cyclic $q$-isopairs on $H^2(\T)\otimes \C^k$ (where as usual
  $\Phi$ and $\Psi$ are $k\times k$ matrix valued rational inner
  functions on $\D$).  Then there exists a matrix valued rational inner
  function $F$ such that
\be \label{conjeq}
\Phi(z) = F(z) \Psi(z) F(z)^{-1}
\ee
\end{cor}
Note that in general an expression like $F \Psi F^{-1}$ need not
be holomorphic in the disk.  

\bp By Theorem \ref{thmc1}, $(M_z,M_\Phi)$ has a restriction to a
finite codimension invariant subspace $\mathcal{K}$ that is unitarily
equivalent to $(M_z,M_\Psi)$.  Since $\mathcal{K}$ is shift invariant
and of finite codimension, it is of the form $F (H^2\otimes \C^k)$ where
$F$ is a rational matrix valued inner function. However, since
$\mathcal{K}$ is invariant under $M_\Phi$, we see that $\Phi F
(H^2\otimes \C^k) \subset F (H^2\otimes \C^k)$.  This implies $G =
F^{-1} \Phi F$ is holomorphic and also rational and inner. It is not
hard to show that $(M_z,M_\Phi)$ on $F(H^2\otimes \C^k)$ is unitarily
equivalent to the pair $(M_z, M_G)$ on $H^2\otimes \C^k$ (\ie the map
$Fg \mapsto g$ is the required Hilbert space isomorphism that
intertwines the operators).  Therefore, $(M_z,M_G)$ and $(M_z,
M_\Psi)$ are unitarily equivalent.  This can only occur if $G$ and
$\Psi$ are unitarily equivalent.  This proves \eqref{conjeq} after
replacing $F$ with an appropriate unitary multiple.
\ep

\section{Convex hulls}
\label{sece}

The operator-theoretic ideas of Sections~\ref{secb} and \ref{secc} allow us to
prove  a result in function theory, Theorem~\ref{thme1} below.
E.L.~Stout has proved a similar result for irreducible analytic subvarieties (private
communication).

The central issue is what one can say about the intersection of an algebraic set $A$ with
the two-torus in $\C^2$. One immediate distinction is whether $A \cap \T^2$ is large in the sense
that no polynomial can vanish on $A \cap \T^2$ without vanishing identically on $A$; if this holds we 
call the set $A$ {\em toral}. (In two dimensions, as we are here, this just means that
$A \cap \T^2$ is infinite, though the definition makes sense in higher dimensions.)
However, the two curves $$
A_1 \ = \ \{ (z,w) \, : \, z=w \}
$$
and $$
A_2 \ = \ \{ (z,w) \, : \, zw = 1 \}
$$
are both toral, yet $X_1 = A_1 \cap \T^2$ and $X_2 = A_2 \cap \T^2$ are
qualitatively different. The first bounds an analytic disk in $\D^2$; the second does not.
Theorem~\ref{thme1} says that one way to understand this is to observe that $A_1$ is a distinguished variety, and
$A_2$ is not.

We wish to exclude curves that contain horizontal or vertical planes
(\ie zero sets of polynomials $z-\z_1$ or $w-\z_2$). If $\z_r$ is
unimodular, such a zero set fills a disk in the boundary of $\D^2$,
and the polynomial hull of the intersection of this disk with the
torus is the closed disk. So, we exclude linear factors to make the
statement of the theorem concise; however one could drop this
restriction and conclude that the set $X$ is not polynomially convex
if and only if $q$ has a factor that is either inner toral or of the
form $(z- e^{i \theta})$ or $(w- e^{i \theta})$.
 
Recall that a polynomial $q$ is inner toral if
$
Z_q \ \subseteq \ \D^2 \cup \T^2 \cup \E^2 .
$

\bt
\label{thme1}
Let $q$ be a polynomial in two variables with no linear factors. Then $Y = Z_q \cap \T^2$ is polynomially 
convex if and only if $q$ has no inner toral factor.
\et

\bp
(i) First assume that $Y$ is not polynomially convex, so there is some point $\z$
in the polynomial hull of $Y$ that is not in $Y$.
As every point of $\T^2$ is a peak point for $A(\D^2)$, we cannot have
$\z$ in $\T^2$.
There exists a complex measure $\l$ supported on $Y$ so that
\be
\label{eqe1}
p(\z) \= \int_Y p \, d \l
\ee
for all polynomials $p$. Let $d \mu = |d \l|$  
be the total variation of $\l$.
Then by (\ref{eqe1}), for every polynomial $p$ we have
\be
\label{eqe2}
|p(\z)| \ \leq \ C\, \left[ \int_Y |p|^2 d \mu \right]^{1/2} .
\ee
(A point $\z$ satisfying inequality (\ref{eqe2}) for all polynomials $p$ is called
a {\em bounded point evaluation for $\ptm$.})

Claim:  $\z \inn \D^2$.

Else, some component, say $\z_1$, is unimodular. Applying (\ref{eqe2}) to 
polynomials of the form
$$
p(z,w) \= \left( \frac{z+\z_1}{2} \right)^n r(w)
$$
and letting $n$ tend to infinity, we would get
\be
\label{eqe3}
|r(\z_2)| \ \leq \ C\, \left[ \int_{Y \cap \{ z = \z_1 \} } |r|^2 d \mu \right]^{1/2} .
\ee
As $\z_2$ is in $\D$, (\ref{eqe3}) 
asserts that the measure $\mu |_{Y \cap \{ z = \z_1 \} }$ is a measure on the circle
that has a bounded point evaluation inside the disk.
By \szs theorem \cite{sze20}, this means
that $Y \cap  \{ z = \z_1 \} $
must be the whole circle $\z_1 \times \T$, and so $q$ must have
$(z-\z_1)$ as a factor, contrary to assumption. 
\oec

Let $V_1$ be $(M_z,M_w)$ on $\ptm$, and let $V$ be the pure part of
$V_1$. By Lemma~\ref{leme1}, $V$ is non-zero.
So, by Theorem~\ref{thmb4}, there is some square-free inner toral polynomial $p$ such that
$p(V) = 0$, and $p = 0 \ \mu_a$-a.e.
As both $p$ and $q$ vanish on the support of $\mu_a$, which is an infinite set,
they must share a common factor, which is an inner toral factor of $q$. 

(ii) Suppose that $q$ has a factor $p$ that is inner toral.
Then $Z_p \cap \D^2$, which is non-empty by Theorem~\ref{thmb1}, is contained in 
the polynomial hull of $Z_p \cap \T^2$, and hence of $Y$.
\ep

\bl
\label{leme1} Suppose $\mu$ is a measure on $Z_q \cap \T^2$, and
$\ptm$ has a bounded point evaluation $\z$ in $\D^2$. Then the isopair
$(M_z,M_w)$ on $\ptm$ has a non-zero pure part, and this is unitarily
equivalent to $(M_z,M_w)$ on $P^2(\mu_a)$, where $\mu_a$ is the part
of $\mu$ that is absolutely continuous with respect to arc-length on
$\T^2$.  
\el 

\bp Using the notation of Section~\ref{secc}, every point in
$h^{-1}(\z)$ (which can be more than one point if $\z$ is a point of
multiplicity of $q$) is a bounded point evaluation for $\athn$, and
therefore for $\atn$.  
Just as in the 
proof of Lemma~\ref{lemcd},
the Kolmogorov-Krein theorem says that 
$$
A^2(\nu) \= A^2(\nu_a) \oplus L^2(\nu_s),
$$
where 
$\nu_a$ and $\nu_s$ are, respectively, the absolutely continuous and singular
parts of $\nu$
with respect to
harmonic measure $\omega$.  
Therefore,
every
point in $S$ is a bounded point evaluation for $A^2(\nu_a)$.
Pushing back down to $\Omega=Z_q\cap\T^2$
again, we find that every point is a bounded point evaluation for
$P^2(\mu_a)$ where $\mu_a = h_* (\nu_a)$.

To see that 
$(M_z,M_w)$ on $P^2(\mu_a)$
is pure, 
assume that one of the isometries, $M_z$ say, has a unitary part.
Then there is some function $f$ of norm one in $P^2(\mu_a)$ such that
\be
\label{eqf4}
\| M_z^{*n} f \| \= \| P \bar z^n f \| \= \| f \|
\ee
for all $n$, 
where $P$ is the projection from  
$L^2(\mu_a)$
onto 
$P^2(\mu_a)$.
From (\ref{eqf4}) we get that 
\be
\label{eqf5}
M_z^{*n} f  \=  P \bar z^n f  \=  \bar z^n f \quad \forall n. 
\ee
Choose some bounded point evaluation $\z = (\z_1,\z_2)$ in $\D^2$ of 
$P^2(\mu_a)$ such that $f(\z) \neq 0$.
(Such a point must exist, for otherwise $f\circ h$ would be in 
$A^2_h(\nu_a)$ and vanish at every point of $S$, and so by \cite{ahesar1} again
would be identically zero.)
Let $k_\z$ be the kernel function at $\z$, \ie the unique function in $P^2(\mu_a)$ satisfying
$$
p(\z) \= \ip{p}{k_\z} \qquad \forall \ {\rm polynomials\ } p .
$$
Then
\beq
f(\z) &\=& \ip{ z^n \bar z^n f}{k_\z} \\
&\=& \ip{ \bar z^n f}{P \bar z^n k_\z} 
\eeq
Therefore
\[
\| P \bar z^n k_\z \| \ \geq |f(\z)| \qquad \forall n .
\]
But 
\[
 P \bar z^n k_\z  \= \bar \z_1^n k_\z ,
\]
and this must tend to zero as $n$ goes to infinity.
\ep

\section{Non-cyclic algebraic isopairs}
\label{secf}

We do not understand algebraic isopairs that are not nearly cyclic.

Let us say that an isopair $V$ is {\em essentially $k$-cyclic} if there are 
$k$ vectors $u_1,\dots, u_k$ so that
$$
\vee \{ p_1(V) u_1, \dots, p_k(V) u_k \ : \ p_1,\dots,p_k \inn \C[z,w] \}
$$
is of finite codimension, and if no set of $k-1$ vectors suffices.

\begin{question}
\label{qf1}
Suppose $V$ and $V'$ are both essentially $k$-cyclic isopairs with the same minimal polynomial.
Are they nearly unitarily equivalent?
\end{question}

\begin{question}
Are all essentially $k$-cyclic algebraic isopairs nearly equivalent to
a direct sum of $k$ cyclic algebraic isopairs?
\end{question}

\bibliography{references}

\end{document}